\documentclass[12pt]{article}
\usepackage{amssymb,amsmath}
\usepackage[dvips]{graphicx}
\usepackage{hyperref}
\usepackage{ulem}

\begin{document}
%\begin{center}
\title{\LARGE\bf A rapid and highly accurate approximation for the error function of complex argument}

%\bigskip
\author{
\normalsize\bf S. M. Abrarov\footnote{\scriptsize{Dept. Earth and Space Science and Engineering, York University, Toronto, Canada, M3J 1P3.}}\, and B. M. Quine$^{*}$\footnote{\scriptsize{Dept. Physics and Astronomy, York University, Toronto, Canada, M3J 1P3.}}}

\date{August 14, 2013}
\maketitle
%\vspace{1cm}%\bigskip

\begin{abstract}
We present efficient approximation of the error function obtained by Fourier expansion of the exponential function $\exp \left[ { - {\left(t - 2 \sigma\right)^2}/4} \right]$. The error analysis reveals that it is highly accurate and can generate numbers that match up to the last decimal digits with reference values. Due to simple representation the proposed error function approximation can be utilized in a rapid algorithm.
\vspace{0.25cm}
\vspace{0.25cm}
\end{abstract}

\section {Complex error function approximation}
Consider the complex error function also known as the Faddeeva function or the Kramp function that is defined as \cite{Faddeeva1961, Schreier1992}
\small
\begin{equation}\label{eq_1}
\begin{aligned}
w\left( z \right) &= {e^{ - {z^2}}}\left( {1 + \frac{{2i}}{\sqrt{\pi} }\int\limits_0^z {{e^{{t^2}}}dt} } \right) \\ 
&= {e^{ - {z^2}}}{\text{erfc}}\left( { - iz} \right), \\ 
\end{aligned}
\end{equation}
\normalsize
where $z = x + iy$ is a complex argument. It is shown that the real and imaginary parts of the complex error function $w\left( {x,y} \right) = K\left( {x,y} \right) + iL\left( {x,y} \right)$ can be expressed as \cite{Srivastava1992}
\begin{equation}\label{eq_2}
K\left( {x,y} \right) = \frac{1}{{\sqrt \pi  }}\int\limits_0^\infty  {\exp \left( { - {t^2}/4} \right)\exp \left( { - yt} \right)\cos \left( {xt} \right)dt}, \quad\quad y > 0
\end{equation}
and
\begin{equation}\label{eq_3}
L\left( {x,y} \right) = \frac{1}{{\sqrt \pi  }}\int\limits_0^\infty  {\exp \left( { - {t^2}/4} \right)\exp \left( { - yt} \right)\sin \left( {xt} \right)dt}, \quad\quad y > 0,
\end{equation}
respectively. The real part \eqref{eq_2} of the complex error function is known as the Voigt function. Combining the real \eqref{eq_2} and imaginary \eqref{eq_3} parts together results in
\begin{equation}\label{eq_4}
w\left( {x,y} \right) = \frac{1}{{\sqrt \pi  }}\int\limits_0^\infty  {\exp \left( { - {t^2}/4} \right)\exp \left( { - yt} \right)\exp \left( {ixt} \right)dt}.
\end{equation}
It is useful to introduce a real constant $\sigma  > 0$ and rewrite the complex error function \eqref{eq_4} as
\begin{equation*}
w\left( {x,y} \right) = \frac{1}{{\sqrt \pi  }}\int\limits_0^\infty  {{e^{ - {t^2}/4}}{e^{ - \left( {y + \sigma } \right)t}}{e^{\sigma t}}{e^{ixt}}dt}
\end{equation*}
or
\begin{equation}\label{eq_5}
w\left( {x,y} \right) = \frac{{{e^{{\sigma ^2}}}}}{{\sqrt \pi  }}\int\limits_0^\infty  {{e^{ - {{\left( {t - 2\sigma } \right)}^2}/4}}{e^{ - \left( {y + \sigma } \right)t}}{e^{ixt}}dt.}
\end{equation}

We have shown previously that exponential function can be approximated in form \cite{Abrarov2010}
\begin{equation}\label{eq_6}
{e^{ - {t^2}/4}} \approx  - \frac{{{a_0}}}{2} + \sum\limits_{n = 0}^N {{a_n}\cos \left( {\frac{{n\pi }}{{{\tau _m}}}t} \right)}, \quad\quad - {\tau _m} \leqslant t \leqslant {\tau _m},
\end{equation}
where the Fourier expansion coefficients are
\begin{equation*}
{a_n} \approx \frac{{2\sqrt \pi  }}{{{\tau _m}}}\exp \left( { - \frac{{{n^2}{\pi ^2}}}{{\tau _m^2}}} \right)
\end{equation*}
and ${\tau _m}$ is the margin value \cite{Abrarov2010, Geetha2010}. Consequently, from approximation \eqref{eq_6} it follows that
\begin{equation}\label{eq_7}
{e^{ - {{\left( {t - 2\sigma } \right)}^2}/4}} \approx  - \frac{{{a_0}}}{2} + \sum\limits_{n = 0}^N {{a_n}\cos \left( {\frac{{n\pi }}{{{\tau _m}}}\left( {t - 2\sigma } \right)} \right)}.
\end{equation}
Geometrically, the peak of the exponential function ${e^{ - {{\left( {t - 2\sigma } \right)}^2}/4}}$ is shifted towards right with respect to origin by value equal to $2\sigma $. Therefore, further we will refer to the value $\sigma $ as the shift constant.

Substituting approximation \eqref{eq_7} into complex error function approximation \eqref{eq_5} leads to the following complex error function approximation
\small
\begin{equation}\label{eq_8}
\begin{aligned}
w\left( {x,y} \right) &= \frac{{{e^{{\sigma ^2}}}}}{{\sqrt \pi  }}\left( {{a_0}\frac{{ix - \left( {y + \sigma } \right)}}{{2{{\left( {x + i\left( {y + \sigma } \right)} \right)}^2}}}} \right. \\ 
&+ \left. {\sum\limits_{n = 1}^N {{a_n}\frac{{\cos \left( {2n\pi \sigma /{\tau _m}} \right)\left( {\left( {y + \sigma } \right) - ix} \right) + n\pi \sin \left( {2n\pi \sigma /{\tau _m}} \right)/{\tau _m}}}{{{{\left( {n\pi /{\tau _m}} \right)}^2} - {{\left( {x + i\left( {y + \sigma } \right)} \right)}^2}}}} } \right).
\end{aligned}
\end{equation}
\normalsize
Approximation \eqref{eq_8} of the complex error function can be conveniently rewritten in form
\begin{equation}\label{eq_9}
w\left( z \right) \approx \frac{{{e^{{\sigma ^2}}}}}{{{\tau _m}\left( {\sigma  - iz} \right)}} + \sum\limits_{n = 1}^N {\frac{{{A_n}\left( {\sigma  - iz} \right) + {B_n}}}{{{n^2}{\pi ^2} + \tau _m^2{{\left( {\sigma  - iz} \right)}^2}}}}, \quad\quad \operatorname{Im} \left[ z \right] > 0,
\end{equation}
where
$$
{A_n} = 2{\tau _m}{e^{{\sigma ^2} - {n^2}{\pi ^2}/\tau _m^2}}\cos \left( {2n\pi \sigma /{\tau _m}} \right)
$$
and
$$
{B_n} = 2n\pi {e^{{\sigma ^2} - {n^2}{\pi ^2}/\tau _m^2}}\sin \left( {2n\pi \sigma /{\tau _m}} \right)
$$
are the coefficients that are independent of $x$ and $y$ parameters.

Another useful approximation of the complex error function is given by \href{http://www.sciencedirect.com/science/article/pii/S0096300311009179}{[6]} (see also \href{http://arxiv.org/pdf/1205.1768.pdf}{[7]})
\footnotesize
\begin{equation}\label{eq_10}
\begin{aligned}
w\left( z \right) &\approx \frac{i}{{2\sqrt \pi  }}\left[ {\sum\limits_{n = 0}^N {{a_n}{\tau _m}\left( {\frac{{1 - {e^{i\,\,\left( {n\pi  + {\tau _m}z} \right)}}}}{{n\,\pi  + {\tau _m}z}} - \frac{{1 - {e^{i\,\,\left( { - n\pi  + {\tau _m}z} \right)}}}}{{n\,\pi  - {\tau _m}z}}} \right)}  - {a_0}\frac{{1 - {e^{i{\tau _m}z}}}}{z}} \right] \\ 
& = i\frac{{1 - {e^{i{\tau _m}z}}}}{{{\tau _m}z}} + i\frac{{\tau _m^2z}}{{\sqrt \pi  }}\sum\limits_{n = 1}^N {{a_n}\frac{{{{\left( { - 1} \right)}^n}{e^{i{\tau _m}z}} - 1}}{{{n^2}\,{\pi ^2} - \tau _m^2{z^2}}},} \quad\quad \operatorname{Im} \left[ z \right] > 0.
\end{aligned}
\end{equation}
\normalsize
In contrast to equation \eqref{eq_10}, the complex error function approximation \eqref{eq_9} is based on superposition of rational functions. Therefore, it is advantageous for more rapid computation.

It should be noted that due to rapid convergence of equations \eqref{eq_6} and \eqref{eq_7}, both complex error function approximations \eqref{eq_9} and \eqref{eq_10} provide essentially higher accuracy than the complex error function approximation proposed earlier by Weideman (see equation (13) in \cite{Weideman1994} and also \cite{Abrarov2011} for comparison). This tendency becomes particularly evident at $y < {10^{ - 3}}$.

\section{Proposed error function approximation}

Let us return to equation \eqref{eq_1} and rewrite it as
\begin{equation}\label{eq_11}
{\text{erfc}}\left( z \right) = {e^{ - {z^2}}}w\left( {iz} \right)
\end{equation}
Consequently, substituting the complex error function approximation \eqref{eq_9} into identity \eqref{eq_11} yields a simple approximation of the complementary error function
\small
\begin{equation}\label{eq_12}
{\text{erfc}}\left( z \right) \approx {e^{ - {z^2}}}\left( {\frac{{{e^{{\sigma ^2}}}}}{{{\tau _m}\left( {\sigma  + z} \right)}} + \sum\limits_{n = 1}^N {\frac{{{A_n}\left( {\sigma  + z} \right) + {B_n}}}{{{n^2}{\pi ^2} + \tau _m^2{{\left( {\sigma  + z} \right)}^2}}}} } \right), \quad\quad \operatorname{Re} \left[ z \right] > 0,
\end{equation}
\normalsize
The corresponding approximation of the error function can be obtained by using the following equation
\begin{equation}\label{eq_13}
{\text{erf}}\left( z \right) = 1 - {\text{erfc}}\left( z \right).
\end{equation}
This leads to
\small
\begin{equation}\label{eq_14}
{\text{erf}}\left( z \right) \approx 1 - {e^{ - {z^2}}}\left( {\frac{{{e^{{\sigma ^2}}}}}{{{\tau _m}\left( {\sigma  + z} \right)}} + \sum\limits_{n = 1}^N {\frac{{{A_n}\left( {\sigma  + z} \right) + {B_n}}}{{{n^2}{\pi ^2} + \tau _m^2{{\left( {\sigma  + z} \right)}^2}}}} } \right), \quad\quad \operatorname{Re} \left[ z \right] > 0.
\end{equation}
\normalsize
 
As all terms inside the brackets of equation \eqref{eq_14} are complex error function approximation $w\left( {iz} \right)$ based on superposition of rational functions, the approximation \eqref{eq_14} is effective for implementation in a rapid algorithm.

\section{Error analysis}

In order to perform error analysis it is convenient to define the relative errors for the real and imaginary parts as
\begin{equation*}
{\Delta _{\operatorname{Re} }} = \left| {\frac{{\operatorname{Re} \left[ {{\text{erf}}\left( z \right)} \right] - \operatorname{Re} \left[ {{\text{er}}{{\text{f}}_{ref.}}\left( z \right)} \right]}}{{\operatorname{Re} \left[ {{\text{er}}{{\text{f}}_{ref.}}\left( z \right)} \right]}}} \right|
\end{equation*}
and
\begin{equation*}
{\Delta _{\operatorname{Im} }} = \left| {\frac{{\operatorname{Im} \left[ {{\text{erf}}\left( z \right)} \right] - \operatorname{Im} \left[ {{\text{er}}{{\text{f}}_{ref.}}\left( z \right)} \right]}}{{\operatorname{Im} \left[ {{\text{er}}{{\text{f}}_{ref.}}\left( z \right)} \right]}}} \right|, 
\end{equation*}
respectively, where ${\text{er}}{{\text{f}}_{ref.}}\left( z \right)$ is the reference. 

To obtain the reference we can substitute the complex error function approximation \eqref{eq_10} into identity \eqref{eq_11}. This results to complimentary error function approximation
\small
\begin{equation}\label{eq_15}
{\text{erfc}}\left( z \right) \approx {e^{ - {z^2}}}\left( {\frac{{1 - {e^{ - {\tau _m}z}}}}{{{\tau _m}z}} + \frac{{\tau _m^2z}}{{\sqrt \pi  }}\sum\limits_{n = 1}^N {{a_n}\frac{{1 - {{\left( { - 1} \right)}^n}{e^{ - {\tau _m}z}}}}{{{n^2}\,{\pi ^2} + \tau _m^2{z^2}}}} } \right), \quad\quad \operatorname{Re} \left[ z \right] > 0.
\end{equation}
\normalsize
Consequently, from equations \eqref{eq_13} and \eqref{eq_15} we obtain highly accurate approximation of the error function suitable to generate the reference values
\small
\begin{equation}\label{eq_16}
{\text{erf}}\left( z \right) \approx 1 - {e^{ - {z^2}}}\left( {\frac{{1 - {e^{ - {\tau _m}z}}}}{{{\tau _m}z}} + \frac{{\tau _m^2z}}{{\sqrt \pi  }}\sum\limits_{n = 1}^N {{a_n}\frac{{1 - {{\left( { - 1} \right)}^n}{e^{ - {\tau _m}z}}}}{{{n^2}\,{\pi ^2} + \tau _m^2{z^2}}}} } \right), \quad\, \operatorname{Re} \left[ z \right] > 0.
\end{equation}
\normalsize

Tables \ref{Table_1} and \ref{Table_2} show the error function values from approximations \eqref{eq_14}, \eqref{eq_16} and relative error for their real and imaginary parts. In both approximations \eqref{eq_14} and \eqref{eq_16} the parameters taken for computation are $N = 23$ and ${\tau _m} = 12$. In approximation \eqref{eq_14} the shift constant is chosen to be $\sigma  = 2$. As we can see from second and third columns in these tables, the numbers generated by approximations \eqref{eq_14} and \eqref{eq_16} can match up to the last decimal digit. The negligibly small relative errors for the real and imaginary parts shown in the last columns of the Tables \ref{Table_1} and \ref{Table_2} confirm high accuracy of the proposed error function approximation \eqref{eq_14}.

%\newpage
\begin{table}[ht]
\caption{\small{Real parts of the error function approximations \eqref{eq_14} and \eqref{eq_16}.\vspace{0.15cm}}\label{Table_1}}
\scriptsize
\centering
\begin{tabular}{c c c c c}
\hline\hline
\bfseries{Parameter x} & \bfseries{Parameter y} & \bfseries{Approximation \eqref{eq_14}} & \bfseries{Approximation \eqref{eq_16}} & \bfseries{$\Delta _{\operatorname{Re}}$ } \\ [0.5ex]
\hline
10 & 10 & 9.616493742724747E-1 & 9.616493742724749E-1 & 2.3090E-16\\
10 & 5 & 1.000000000000000E0 & 1.000000000000000E0 & 0\\
5 & 5 & 9.303796037430947E-1 & 9.303796037430951E-1 & 4.7732E-16\\
5 & 1 & 1.000000000002960E0 & 1.000000000002960E0 & 0\\
1 & 1 & 1.316151281697949E0 & 1.316151281697948E0 & 6.7483E-16\\
1 & 0.5 & 9.507097283189570E-1 & 9.507097283189572E-1 & 2.3356E-16\\
0.5 & 0.5 & 6.426129148548198E-1 & 6.426129148548205E-1 & 1.2094E-15\\
0.5 & 0.1 & 5.249121488205361E-1 & 5.249121488205371E-1 & 1.9036E-15\\
0.1 & 0.1 & 1.135856345618654E-1 & 1.135856345618664E-1 & 8.7969E-15\\
0.1 & 0.05 & 1.127425509896926E-1 & 1.127425509896922E-1 & 2.9542E-15\\
0.05 & 0.05 & 5.651284873688534E-2 & 5.651284873688744E-2 & 3.7326E-14\\
0.05 & 0.01 & 5.637760588665819E-2 & 5.637760588665963E-2 & 2.5600E-14\\
0.01 & 0.01 & 1.128454387859423E-2 & 1.128454387859545E-2 & 1.0822E-13\\
0.01 & 0.005 & 1.128369762595771E-2 & 1.128369762595882E-2 & 9.8392E-14\\
0.005 & 0.005 & 5.641989865663000E-3 & 5.641989865664443E-3 & 2.5581E-13\\
0.005 & 0.001 & 5.641854461787554E-3 & 5.641854461787776E-3 & 3.9357E-14\\
0.001 & 0.001 & 1.128379919345890E-3 & 1.128379919343114E-3 & 2.4598E-12\\[1ex]
\hline
\end{tabular}
\normalsize
\end{table}

\newpage
\begin{table}[ht]
\caption{\small{Imaginary parts of the error function approximations \eqref{eq_14} and \eqref{eq_16}.\vspace{0.15cm}}\label{Table_2}}
\scriptsize
\centering
\begin{tabular}{c c c c c}
\hline\hline
\bfseries{Parameter x} & \bfseries{Parameter y} & \bfseries{Approximation \eqref{eq_14}} & \bfseries{Approximation \eqref{eq_16}} & \bfseries{$\Delta _{\operatorname{Im}}$ } \\ [0.5ex]
\hline
10 & 10 & -1.098768460819404E-2 & -1.098768460819399E-2 & 4.2627E-15\\
10 & 5 & -9.495949264558077E-36 & -9.495949264558098E-36 & 2.2517E-15\\
5 & 5 & 3.893619089512146E-2 & 3.893619089512138E-2 & 2.1385E-15\\
5 & 1 & -2.846018382085604E-12 & -2.846018382085594E-12 & 3.5479E-15\\
1 & 1 & 1.904534692378354E-1 & 1.904534692378347E-1 & 3.2062E-15\\
1 & 0.5 & 1.879734672233839E-1 & 1.879734672233833E-1 & 3.2485E-15\\
0.5 & 0.5 & 4.578813944351928E-1 & 4.578813944351922E-1 & 1.4548E-15\\
0.5 & 0.1 & 8.802479434588868E-2 & 8.802479434588850E-2 & 2.0496E-15\\
0.1 & 0.1 & 1.120811719910652E-1 & 1.120811719910650E-1 & 1.9811E-15\\
0.1 & 0.05 & 5.590323090214489E-2 & 5.590323090214489E-2 & 0\\
0.05 & 0.05 & 5.632478587819852E-2 & 5.632478587819856E-2 & 6.1597E-16\\
0.05 & 0.01 & 1.125599074671486E-2 & 1.125599074671481E-2 & 5.2399E-15\\
0.01 & 0.01 & 1.128303937304405E-2 & 1.128303937304429E-2 & 2.1832E-14\\
0.01 & 0.005 & 5.641378676150062E-3 & 5.641378676150111E-3 & 8.7637E-15\\
0.005 & 0.005 & 5.641801802469853E-3 & 5.641801802469647E-3 & 3.6436E-14\\
0.005 & 0.001 & 1.128351334067245E-3 & 1.128351334067475E-3 & 2.0351E-13\\
0.001 & 0.001 & 1.128378414842284E-3 & 1.128378414846887E-3 & 4.0794E-12\\[1ex]
\hline
\end{tabular}
\normalsize
\end{table}

\section{Extension to the entire complex plane}

Approximations \eqref{eq_12}, \eqref{eq_14}, \eqref{eq_15} and \eqref{eq_16} are valid only when the parameter $x = \operatorname{Re} \left[ z \right]$  is positive. Despite of it, these approximations are applicable not only in the first $x > 0,y > 0$ and fourth $x > 0,y < 0$ quadrants, but can also be extended to the second $x < 0,y > 0$ and third $x < 0,y < 0$ quadrants. Specifically, from the property of the complex error function \cite{Armstrong1967, Zaghloul2011}
$$
w\left( { - z} \right) = 2{e^{ - {z^2}}} - w\left( z \right)
$$ 
and equation \eqref{eq_1} it follows that
$$
{\text{erfc}}\left( { - z} \right) = 2 - {\text{erfc}}\left( z \right),
$$
or
\begin{equation}\label{eq_17}
{\text{erfc}}\left( { - x \mp iy} \right) = 2 - {\text{erfc}}\left( {x \pm iy} \right)
\end{equation}
Finally, from equation \eqref{eq_13} and \eqref{eq_17} we obtain
\begin{equation}\label{eq_18}
- {\text{erf}}\left( { - x \mp iy} \right) = {\text{erf}}\left( {x \pm iy} \right).
\end{equation}
Considering equations \eqref{eq_17} and \eqref{eq_18} we can notice that the parameter $x$ on the right side is positive. This signifies that the approximations \eqref{eq_12} and \eqref{eq_14} of the complimentary error function can be applied for the entire coverage of the complex plane according to equation \eqref{eq_17}. Similarly, the approximations (15) and (16) can be employed for the entire coverage of the complex plane according to equation \eqref{eq_18}.

\section{Conclusion}

We present an error function approximation obtained by Fourier expansion of the exponential function $\exp \left[ { - {{\left( {t - 2\sigma } \right)}^2}/4} \right]$ for high-accuracy computation. The error analysis shows that the proposed approximation \eqref{eq_14} generates numbers matching up to the last decimal digits with reference values. Due to simple representation, the proposed error function approximation is efficient for implementation in a rapid algorithm.

\section*{Acknowledgements}
This work is supported by the National Research Council of Canada, Thoth Technology Inc., and York University. The authors wish to thank Prof. Ian McDade and Dr. Brian Solheim for discussions and constructive suggestions.

%\bigskip
%\newpage


\begin{thebibliography}{9}
%\small
\bibitem{Faddeeva1961}
V. N. Faddeyeva and N. M. Terent’ev, Tables of the probability integral for complex argument, Pergamon Press, Oxford, 1961.

\bibitem{Schreier1992}
F. Schreier, The Voigt and complex error function: A comparison of computational methods, J. Quant. Spectrosc. Radiat. Transfer, 48 (1992) 743-762.

\bibitem{Srivastava1992}
H. M. Srivastava and M. P. Chen, Some unified presentations of the Voigt functions, Astrophys. Space Sci., 192 (1992) 63-74.

\bibitem{Abrarov2010}
S. M. Abrarov, B. M. Quine and R. K. Jagpal, Rapidly convergent series for high-accuracy calculation of the Voigt function, J. Quant. Spectrosc. Radiat. Transfer, 111 (2010) 372-375.

\bibitem{Geetha2010}
R. S. Geetha, R. S. Keshavamurthy and R. Harish, Temperature and energy derivatives of Doppler broadening functions by Steffensen’s inequality technique, Ann. Nucl. Energy, 37 (2010) 985-990.

\bibitem{Abrarov2011}
S. M. Abrarov and B. M. Quine, Efficient algorithmic implementation of the Voigt/complex error function based on exponential series approximation, \href{http://www.sciencedirect.com/science/article/pii/S0096300311009179}{Appl. Math. Comput., 218 (2011) 1894-1902.}

\bibitem{Abrarov2012}
S. M. Abrarov and B. M. Quine, On the Fourier expansion method for highly accurate computation of theVoigt/complex error function in a rapid algorithm, \href{http://arxiv.org/pdf/1205.1768.pdf}{arXiv:1205.1768v1}

\bibitem{Weideman1994}
J. A. C. Weideman, Computation of the complex error function, SIAM J. Numer. Anal. 31 (1994) 1497-1518.

\bibitem{Armstrong1967}
B. H. Armstrong, Spectrum line profiles: The Voigt function, J. Quant. Spectrosc. Radiat. Transfer, 7 (1967) 61-88.

\bibitem{Zaghloul2011}
M. R. Zaghloul and A. N. Ali, Algorithm 916: computing the Faddeyeva and Voigt functions, ACM Trans. Math. Software 38 (2011) 15:1-15:22.

\end{thebibliography}
\end{document}